
\def\LaTeX{%
  \let\Begin\begin
  \let\End\end
  \let\salta\relax
  \let\finqui\relax
  \let\futuro\relax}

\def\UK{\def\our{our}\let\sz s}
\def\USA{\def\our{or}\let\sz z}



\LaTeX

\USA


\salta

\documentclass[10pt,a4paper]{article}  
\setlength{\textheight}{205mm}
\setlength{\textwidth}{145mm}
\setlength{\oddsidemargin}{1cm}
\setlength{\topmargin}{12mm}
\parskip2mm


\usepackage{color}
\usepackage{amsmath}
\usepackage{amsthm}
\usepackage{amssymb}
\usepackage[mathcal]{euscript}








\bibliographystyle{plain}


%

\finqui

\def\Beq{\Begin{equation}}
\def\Eeq{\End{equation}}
\def\Bsist{\Begin{eqnarray}}
\def\Esist{\End{eqnarray}}

\def\Bthm{\Begin{theorem}}
\def\Ethm{\End{theorem}}
\def\Blem{\Begin{lemma}}
\def\Elem{\End{lemma}}

\def\Brem{\Begin{remark}\rm}
\def\Erem{\End{remark}}

\def\Bdim{\Begin{proof}}
\def\Edim{\End{proof}}
\let\non\nonumber




\def\step #1 \par{\medskip\noindent{\bf #1.}\quad}


\def\Lip{Lip\-schitz}

\def\aand{\quad\hbox{and}\quad}
\def\where{\quad\hbox{where}\quad}

\def\lhs{left-hand side}
\def\rhs{right-hand side}
\def\sfw{straightforward}
\def\omegalimit{$\omega$-limit}


\def\characteriz{characteri\sz}

\def\generaliz{generali\sz}

\def\regulariz{regulari\sz}

\def\nbh{neighb\our hood}
\def\bhv{behavi\our}


\def\multibold #1{\def\arg{#1}%
  \ifx\arg\pto \let\next\relax
  \else
  \def\next{\expandafter
    \def\csname #1#1#1\endcsname{{\bf #1}}%
    \multibold}%
  \fi \next}

\def\pto{.}

\def\multical #1{\def\arg{#1}%
  \ifx\arg\pto \let\next\relax
  \else
  \def\next{\expandafter
    \def\csname cal#1\endcsname{{\cal #1}}%
    \multical}%
  \fi \next}

\def\multizero #1{\def\arg{#1}%
  \ifx\arg\pto \let\next\relax
  \else
  \def\next{\expandafter
    \def\csname #1z\endcsname{#1_0}%
    \multizero}%
  \fi \next}

\def\multistar #1{\def\arg{#1}%
  \ifx\arg\pto \let\next\relax
  \else
  \def\next{\expandafter
    \def\csname #1star\endcsname{#1^*}%
    \multistar}%
  \fi \next}

\let\hat\widehat

\def\multihat #1{\def\arg{#1}%
  \ifx\arg\pto \let\next\relax
  \else
  \def\next{\expandafter
    \def\csname hat#1\endcsname{\hat #1}%
    \multihat}%
  \fi \next}


\def\multimathop #1 {\def\arg{#1}%
  \ifx\arg\pto \let\next\relax
  \else
  \def\next{\expandafter
    \def\csname #1\endcsname{\mathop{\rm #1}\nolimits}%
    \multimathop}%
  \fi \next}

\multibold
qwertyuiopasdfghjklzxcvbnmQWERTYUIOPASDFGHJKLZXCVBNM.

\multical
QWERTYUIOPASDFGHJKLZXCVBNM.

\multizero
qwertyuiopadfghjklzxcvbnmQWERTYUIOPASDFGHJKLZXCVBNM.  

\multistar
qwertyuiopasdfghjklzxcvbnmQWERTYUIOPASDFGHJKLZXCVBNM.

\multihat
qwertyuiopasdfghjklzxcvbnmQWERTYUIOPASDFGHJKLZXCVBNM.

\multimathop
dist div dom meas sign supp .


\def\accorpa #1#2{\eqref{#1}--\eqref{#2}}
\def\Accorpa #1#2 #3 {\gdef #1{\eqref{#2}--\eqref{#3}}%
  \wlog{}\wlog{\string #1 -> #2 - #3}\wlog{}}


\def\tonde #1{\left(#1\right)}

\def\graffe #1{\mathopen\{#1\mathclose\}}
\def\Graffe #1{\left\{#1\right\}}
\def\<#1>{\mathopen\langle #1\mathclose\rangle}
\def\norma #1{\mathopen \| #1\mathclose \|}

\let\meno\setminus

\def\iot {\int_0^t}

\def\ios {\int_0^s}
\def\intQt{\int_{Q_t}}
\def\intQs{\int_{Q_s}}
\def\intQ{\int_{Q_T}}
\def\intQi{\int_{Q_\infty}}
\def\iO{\int_\Omega}

\def\dt{\partial_t}
\def\dn{\partial_n}
\def\ds{\,ds}
\def\dtau{\,d\tau}

\def\cpto{\,\cdot\,}

\def\checkmmode #1{\relax\ifmmode\hbox{#1}\else{#1}\fi}
\def\aeO{\checkmmode{a.e.\ in~$\Omega$}}
\def\aeQ{\checkmmode{a.e.\ in~$Q_T$}}
\def\aeQi{\checkmmode{a.e.\ in~$Q_\infty$}}

\def\aaO{\checkmmode{for a.a.~$x\in\Omega$}}
\def\aaQ{\checkmmode{for a.a.~$(x,t)\in Q_T$}}
\def\aaQi{\checkmmode{for a.a.~$(x,t)\in Q_\infty$}}

\def\aat{\checkmmode{for a.a.~$t\in(0,T)$}}

\def\limn{\lim_{n\to\infty}}


\def\erre{{\mathbb{R}}}




\def\genspazio #1#2#3#4#5{#1^{#2}(#5,#4;#3)}
\def\spazio #1#2#3{\genspazio {#1}{#2}{#3}T0}

\def\spazioi #1#2#3{\genspazio {#1}{#2}{#3}\infty 0}
\def\L {\spazio L}
\def\H {\spazio H}
\def\W {\spazio W}

\def\LL {\spazioi L}

\def\C #1#2{C^{#1}([0,T];#2)}
\def\Ct #1#2{C^{#1}([0,t];#2)}
\def\Cs #1#2{C^{#1}([0,s];#2)}


\def\Lx #1{L^{#1}(\Omega)}
\def\Hx #1{H^{#1}(\Omega)}
\def\Wx #1{W^{#1}(\Omega)}
\def\Luno{\Lx 1}
\def\Ldue{\Lx 2}
\def\Linfty{\Lx\infty}
\def\Huno{\Hx 1}
\def\Hdue{\Hx 2}

\def\Lloc{L^1_{\rm loc}[0,+\infty)}


\def\LQs #1{L^{#1}(Q_s)}
\def\LQ #1{L^{#1}(Q_T)}
\def\LQi #1{L^{#1}(Q_\infty)}


\let\theta\vartheta
\let\eps\varepsilon
\let\phi\varphi

\let\TeXchi\chi                         
\newbox\chibox
\setbox0 \hbox{\mathsurround0pt $\TeXchi$}
\setbox\chibox \hbox{\raise\dp0 \box 0 }
\def\chi{\copy\chibox}


\def\normaV #1{\norma{#1}_V}
\def\normaH #1{\norma{#1}_H}

\let\u\rho
\def\uz{\u_0}
\def\umin{{\u_*}}
\def\umax{{\u^*}}
\def\squ{\sqrt\u}

\def\xiz{\xi_0}
\def\ximax{{\xi^*}}
\def\ximaxz{{\xi^*_0}}
\def\sqx{\sqrt\xi}

\def\sqxz{\sqrt{\xiz}}

\def\DPhi{D(\Phi)}

\def\emmep{{\bar\emme}^+}
\def\emmem{{\bar\emme}^-}
\let\emme\sigma

\def\aeps{\alpha_\eps}
\def\ueps{\u_\eps}
\def\zetaeps{\zeta_\eps}

\def\fu{f_1}
\def\fd{f_2}

\def\uomega{\u^\infty}
\def\ui{\u_\infty}
\def\phii{\phi_\infty}
\def\zetai{\zeta_\infty}

\def\tn{t_n}
\def\un{\u_n}
\def\phin{\phi_n}
\def\zetan{\zeta_n}


\def\cd{c_\delta}
\def\Md{M_2}
\def\Rp{R_p}
\def\Rinfty{R_\infty}
\def\Cstar{C_*}

\DeclareMathAlphabet{\mathbf}{OT1}{cmr}{bx}{it}

\newcommand{\0}{\mathbf 0}
\newcommand{\hb}{\mathbf{h}}
 \font\mba=cmmib10 scaled
\magstephalf

\def\csib{\hbox{\mba {\char 24}}}




\Begin{document}


\title{Global solution and long-time \bhv\\
  for a problem of phase segregation of the Allen-Cahn type} 

\author{%
\normalsize Pierluigi Colli$^{(1)}$\\
{\footnotesize e-mail: {\tt pierluigi.colli@unipv.it}}
\\
\normalsize Gianni Gilardi$^{(1)}$\\
{\footnotesize e-mail: {\tt gianni.gilardi@unipv.it}}
\\
\normalsize Paolo Podio-Guidugli$^{(2)}$\\
{\footnotesize e-mail: {\tt ppg@uniroma2.it}}
\\
\normalsize J\"urgen Sprekels$^{(3)}$\\
{\footnotesize e-mail: {\tt sprekels@wias-berlin.de}}\\
\\
\\
$^{(1)}$
{\small Dipartimento di Matematica ``F. Casorati'', Universit\`a di Pavia}\\
{\small via Ferrata 1, 27100 Pavia, Italy}
\\
$^{(2)}$
{\small Dipartimento di Ingegneria Civile, Universit\`a di Roma ``Tor Vergata''}\\
{\small via del Politecnico 1, 00133 Roma, Italy}
\\
$^{(3)}$
{\small 
Weierstra\ss-Institut f\"ur Angewandte Analysis und Stochastik}\\
{\small Mohrenstra\ss e\ 39, 10117 Berlin, Germany}}

\date{}

\maketitle

\Begin{abstract}
In this paper we study a model for phase segregation consisting
in a sistem of a partial and an ordinary differential equation. 
By a careful definition of maximal solution to the latter equation, this system reduces to
an Allen-Cahn equation with a memory term. Global existence 
and uniqueness of a smooth solution are proven and a characterization 
of the \omegalimit\ set is given. 
 
\vskip3mm

\vskip3mm

\End{abstract}

\salta

\pagestyle{myheadings}
\markright{}

\finqui



\section{Problem setting}
The Allen-Cahn equation
\begin{equation}\label{AC}
 \kappa \, \dt\u - \Delta\u + f'(\u) =0
 \end{equation}
is meant to describe evolutionary processes in a two-phase material body, including \emph{phase segregation}: 
$\rho$, with $\rho(x,t)\in[0,1]$, is an \emph{order-parameter} field 
interpreted as the scaled volumetric density of one of the two phases, 
$\kappa>0$~is a \emph{mobility} coefficient, and $f$ denotes a double-well potential confined in~$(0,1)$ and singular at endpoints. 
The derivation of this equation offered by Gurtin in \cite{gurtin} is based on a \textit{balance of contact
and distance microforces}:
\begin{equation}\label{balance}
\div\csib+\pi+\gamma=0
\end{equation}
and on a `purely mechanical'
\textit{dissipation inequality} restricting the free-energy
growth:
\begin{equation}\label{dissipation}
\partial_t\psi\leq w,\quad w:=-\pi\,\partial_t\rho+\csib\cdot\nabla(\partial_t\rho),
\end{equation}
where the distance microforce is split in an internal part $\pi$ and an external part $\gamma$, $\csib$ denotes the \emph{microscopic stress} vector, and $w$ specifies the (distance and contact) internal microworking;\footnote{{\color{black}In \cite{fremond} the balance of microforces is stated under form of a principle of virtual power 
for microscopic motions.}} the Coleman-Noll compatibility of the constitutive choices 
\begin{equation}\label{constitutive}
\pi=\hat\pi(\rho,\nabla\rho,\partial_t\rho), \quad\csib=\hat\csib(\rho,\nabla\rho,\partial_t\rho),\quad \textrm{and}\quad \psi=\hat\psi(\rho,\nabla\rho)=f(\rho)+\frac{1}{2}|\nabla\rho|^2
\end{equation}
with the dissipation inequality (\ref{dissipation}) yields 
\begin{equation}\label{picsi}
\hat\pi(\rho,\nabla\rho,\partial_t\rho)=-f'(\rho)-\hat\kappa(\rho,\nabla\rho,\partial_t\rho)\partial_t\rho,\quad \hat\csib(\rho,\nabla\rho,\partial_t\rho)=\nabla\rho,
\end{equation}
and hence the Allen-Cahn equation (\ref{AC}), for $\hat\kappa(\rho,\nabla\rho,\partial_t\rho)=\kappa$ and $\gamma\equiv 0$. 

One of us proposed in \cite{Podio} a modified version of Gurtin's derivation, where the dissipation inequality (\ref{dissipation}) is dropped and the microforce balance (\ref{balance}) is coupled with the \emph{microenergy balance}
\begin{equation}\label{energy}
\partial_t\varepsilon=e+w,\quad e:=-\div{\bar \hb}+{\bar \sigma},
\end{equation}
and the \emph{microentropy imbalance}
\begin{equation}\label{entropy}
\partial_t\eta\geq -\div\hb+\sigma,\quad \hb:=\mu{\bar \hb},\quad \sigma:=\mu\,{\bar \sigma}.
\end{equation}
The salient new feature of this approach to phase-segregation modeling is that the \emph{microentropy inflow} $(\hb,\sigma)$ is deemed proportional to the \emph{microenergy inflow} $({\bar \hb},{\bar\sigma})$ through the \emph{chemical potential} $\mu$, a positive field; consistently, the free energy is defined to be
\begin{equation}\label{freeenergy}
\psi:=\varepsilon-\mu^{-1}\eta,
\end{equation}
with the chemical potential playing the same role as \emph{coldness} in the deduction of the heat equation.\footnote{{\color{black}Just as absolute temperature can be seen as a macroscopic measure of microscopic \emph{agitation}, its inverse - the coldness - measures microscopic \emph{quiet}; likewise, the chemical potential can be seen as a macroscopic measure of microscopic \emph{organization}.}}  Combination of (\ref{energy})-(\ref{freeenergy}) gives:
\begin{equation}\label{reduced}
\partial_t\psi\leq -\eta_{}(\mu^{-1})\dot{}+\mu^{-1}{\bar \hb}\cdot\nabla\mu-\pi\,\partial_t\rho+\csib\cdot\nabla(\partial_t\rho),
\end{equation}
an inequality that replaces for (\ref{dissipation}) in restricting constitutive choices that can now be more general than those in (\ref{constitutive}). On taking all of the constitutive mappings delivering $\pi,\xi,\eta$, and ${\bar h}$ depending on the list $\rho,\nabla\rho,\partial_t\rho,$ \emph{and} the chemical potential $\mu$, and on choosing
\begin{equation}\label{constitutives}
\psi=\hat\psi(\rho,\nabla\rho,\mu)=-\mu\,\rho+f(\rho)+\frac{1}{2}|\nabla\rho|^2,
\end{equation}
compatibility with (\ref{reduced}) yields
\begin{equation}\label{cn}
\begin{array}{cl}
\hat\pi(\rho,\nabla\rho,\partial_t\rho,\mu)&\displaystyle{\!\!\!=\mu-f'(\rho)-\hat\kappa(\rho,\nabla\rho,\partial_t\rho)\partial_t\rho,\quad \hat\csib(\rho,\nabla\rho,\partial_t\rho,\mu)=\nabla\rho,}\\ [2.mm]\hat\eta(\rho,\nabla\rho,\partial_t\rho,\mu)&\displaystyle{\!\!\!=-\mu^2\rho,\quad \hat{\bar \hb}(\rho,\nabla\rho,\partial_t\rho,\mu)\equiv \0.}
\end{array}
\end{equation}
With the use of (\ref{cn}) and, once again,  of the additional constitutive assumptions that the mobility $\kappa$ is a positive constant and the external distance microforce $\gamma$ is null, the microforce balance (\ref{balance}) and the energy balance (\ref{energy}) become, respectively, 
\begin{equation}  \label{Iprima}
\kappa \, \dt\u - \Delta\u + f'(\u) = \mu
\end{equation}
and
\begin{equation} \label{secondanew}
 \dt(-\mu^2\rho) = \mu\left(\kappa \, (\dt\u)^2 + \bar\emme\right).
\end{equation}

This nonlinear system consists of a parabolic PDE and a first-order-in-time ODE and is to be solved for the order-parameter field $\rho$ and the chemical potential field $\mu$; formally, setting $\mu\equiv 0$ restitutes the standard Allen-Cahn equation (\ref{AC}). We supplement system (\ref{Iprima})-(\ref{secondanew}) with
the homogeneous Neumann condition
\Beq   \label{Ibcic1}
  \dn \u = 0 \quad \hbox{on the body's boundary} 
  \Eeq
(here $\dn$ denotes the outward normal derivative) and with the initial conditions 
\Beq \label{Ibcic2}
   \u|_{t=0} = \uz\quad\hbox{bounded away from} \;0 \,,\quad\mu|_{t=0} = \mu_0\geq 0\,.
\Eeq
Note that, in view of the third of relations (\ref{cn}), the microentropy cannot exceed the level $0$ from below anywhere at any time, and that the corresponding prescribed initial field 
\Beq\label{eta0}
\eta|_{t=0} =\eta_0=-\mu_0^2\rho_0
\Eeq
is nonpositive-valued.
\vskip 6pt
\Brem
\label{enerflusso}
The last of (\ref{cn}) implies that both the energy influx and the entropy influx are everywhere null for all times. This result, that  has a pivotal role in reducing the energy balance to an ODE, is a direct consequence of assuming that the energy-influx mapping $\hat{\bar\hb}$, just as all the other constitutive mappings, be independent of the gradient of the chemical potential. This assumption, that precludes energy and entropy diffusion, does not seem appropriate in the case of higher-order phase segregation models of Cahn-Hilliard type, like the one proposed in \cite{Podio}. 
\Erem


\Accorpa\Ipbcompl Iprima Ibcic2

\section{Solution strategy and summary of contents}
\setcounter{equation}{0}
The aim of our paper is a mathematical investigation of problem \Ipbcompl. The key idea of our strategy is to attack the problem sequentially, the ODE first, then the PDE. 

To do so, we introduce a change of variable that is expedient to give (\ref{secondanew}) plus (\ref{eta0}) the form of a parametric initial-value problem. The change of variable in question is:
%
\Beq\label{defxi}
\xi := -\eta,\quad \xi_0:=-\eta_0,
\Eeq
whence 
\Beq\label{defmu}
\mu= \sqrt{ \xi / \u };
\Eeq
it leads to 
\Beq
\dt\xi  + \frac {\kappa \, (\dt\u)^2 + \bar\emme} \squ \, \sqx
\label{1pier}=0,\quad \xi|_{t=0} =\xi_0,
\Eeq
a Cauchy problem for $\xi(x,\cdot)$ parameterized on the space variable $x$ and on the field $\rho(x,\cdot)$. 
The general form of this problem is discussed in the Appendix. Suffice it to note here that \eqref{1pier} exhibits the 
Peano phenomenon and has infinitely many solutions; among
them, 
we pick a suitably defined \emph{maximal solution} $\xi $ (or $\sqrt{\xi}$, see definition (\ref{seconda})-(\ref{defS})), having the desirable 
property to stay positive as long as is possible. 
Next, we transform \eqref{Iprima} into
\Beq
 \kappa \, \dt\u - \Delta\u + f'(\u) -  \sqx \frac 1\squ = 0,
\label{2pier}
\Eeq
that is, an Allen-Cahn equation for $\rho(x,\cdot)$ with the additional term 
$- \sqrt{ \xi / \u }$; since the factor $\sqrt{ \xi }$ is implicitly defined in terms of $\u$ as the maximal solution of  \eqref{1pier}, 
\eqref{2pier} must be regarded as an \emph{integrodifferential equation}. 
We prove existence, regularity and uniqueness of the solution to (\ref{2pier}) subject to the boundary condition 
(\ref{Ibcic1}) and the initial condition $(\ref{Ibcic2})_1$ by means of a fixed-point
argument, taking advantage of the iterated Contraction Mapping Principle. Crucial to success is that $\dt\u$  be a priori
uniformly  bounded in the space-time domain; we show that this is the case by applying standard regularity arguments
for parabolic equations.

For a detailed discussion of the  
problem transformation sketched here above, and for a precise exposition of our main results,  
we refer the reader to Section~\ref{MainResults}. 
Our well-posedness results are proved in Section~\ref{WELLPOSEDNESS}. 
Section~\ref{LONGTIME} is devoted to an investigation of the long-time behavior of the solution; we prove that $\sqrt{ \xi}$ uniquely converges to 
some function $\varphi_\infty$ and 
that any element $\u_\infty $ of the $\omega-$limit set solves the stationary problem
\Beq
 \kappa \, \dt\u - \Delta\u_\infty + f'(\u_\infty) -  \varphi_\infty \frac 1{\sqrt {\u_\infty}} = 0,
\label{3pier}
\Eeq
supplemented by suitable homogeneous Neumann boundary conditions.


\section{Main results}
\label{MainResults}
\setcounter{equation}{0}

We begin by specifying once and for all the class of data we consider; further assumptions of local importance will be stated when needed. Our problem is formulated over a space-time cylinder
\Beq
  Q_T = \Omega \times [0,T)
  \quad \hbox{with $T\in(0,+\infty)$},
 \Eeq
where $\Omega$~is an open, bounded and connected set of~$\erre^N (N\geq1)$, with a smooth boundary~$\Gamma$ (we use the notation ~$Q_t:= \Omega \times [0,t)$ for every $t\in(0,+\infty]$). As to the \emph{coarse-grain free energy} $f$, we split it as follows:
\Bsist
  && 0 \leq f = \fu + \fd,
  \where
  \fu, \fd : (0,1) \to \erre
  \quad \hbox{are $C^2$-functions,}
  \label{hpfa}
  \\
  && \hbox{$\fu$ is convex}, \quad
  \hbox{$\fd'$ is bounded}, \quad
  \lim_{r\searrow 0} f'(r) = - \infty ,
  \aand
  \lim_{r\nearrow 1} f'(r) = + \infty .
  \qquad
  \label{hpfb}
\Esist
\Accorpa\Hpf hpfa hpfb
We regard $\fd$ as a smooth perturbation of the singular convex  part $\fu$ of $f$, which is well exemplified by
\Beq
  \fu(r) = r \ln r + (1-r) \ln (1-r)
  \quad \hbox{for $r\in (0,1)\,.$}
  \label{logpot}
\Eeq
As to the energy source and the initial data, we assume that
\Beq
  \bar\emme \in \LQ2 , \quad
  \uz , \xiz \in \Linfty , \quad
  0 < \uz < 1 \aand \xiz \geq 0
  \quad \aeO .
  \label{hpdati}
\Eeq
Finally, we recall that the mobility $\kappa$~is a given positive constant.

With this, we take up the forward Cauchy problem~(\ref{1pier}). 
Clearly, $\xi$~must be nonnegative.
We notice that, if we look for a strictly positive~$\xi$
(for given~$\rho>0$ and $\xiz>0$),
the Cauchy problem~(\ref{1pier}) has a unique local solution.
On the contrary, uniqueness is no longer guaranteed
if we allow $\xi$ to be just nonnegative.
On the other hand, every nonnegative local solution can be extended
to a global solution.
Therefore, we select~a (global) solution to problem (\ref{1pier})
according to the following \emph{maximality criterion}
(for a justification, see the Appendix):
\Bsist
  && \sqrt{\xi(x,t)}
  = \sup \Graffe{w(x,t) :\ w\in\calS^*(\bar\emme,\xiz,\u)}
  \quad \hbox{for $(x,t)\in Q_T,$}
  \where
  \label{seconda}
  \\
  && \calS^*(\bar\emme,\xiz,\u) := \Bigl\{
  w\in\W{1,1}\Luno :\ w(0) = \sqxz, \quad w \geq 0 \quad \aeQ, 
  \qquad
  \non
  \\
  && \phantom{\calS := \Bigl\{} \quad 
    \dt w = - \bigl( \kappa \, (\dt\u)^2 + \bar\emme \bigr) / (2\u^{1/2})
    \quad \hbox{a.e.\ where $w>0$}
  \Bigr\} .
  \label{defS}
\Esist
Accordingly, the maximal $\xi$~satisfies:
\Beq   \label{secondaa}
 \sqrt{\xi(x,t)}
  = \sqrt{\xiz(x)} - \iot \astar(x,s) \ds\,,
  \Eeq
 where
\Beq \label{secondab}
 \astar(x,s) :=\left\{
 \begin{array}{l}
 \displaystyle{\frac{\kappa \, |\dt\u(x,s)|^2 + \bar\emme(x,s)}{2 \, \sqrt{\u(x,s)}}}
 \phantom 0
 \quad \hbox{if} \ \ \xi(x,s) > 0, \\
 0
 \hphantom{\displaystyle{\frac{\kappa \, |\dt\u(x,s)|^2 + \bar\emme(x,s)}{2 \, \sqrt{\u(x,s)}}}}
 \quad \hbox{otherwise}.
 \end{array}
 \right.
\Eeq
\Accorpa\Seconda seconda secondab

At this point, we replace $\mu$ by $\sqrt{\xi/\u}$ in~\eqref{Iprima}
 and supplement the equation (\ref{2pier}) we get
with the boundary and initial conditions for $\u$
given by, respectively, \eqref{Ibcic1} and the first of~\eqref{Ibcic2}. 
Of the so-obtained initial/boundary value problem we give the following variational formulation:

\noindent for
\Beq
  V := \Huno
  \aand
  H := \Ldue,
  \label{defspazi}
\Eeq
seek a field $\u$ such that:
\Bsist
  && \u \in \H1H \cap \C0V\,; 
  \label{regua}
  \\
  && \u(0) = \uz,\quad 0 < \u < 1 \quad \aeQ,\quad
  \frac 1\u + \frac 1{1-\u} \in \LQ\infty \,;
  \label{regub}
  \\
  && \kappa \iO \dt\u(t) \, z
  + \iO \nabla\u(t) \cdot \nabla z
  + \iO f'(\u(t)) \, z
  - \iO \bigl( \xi(t)/\u(t) \bigr)^{1/2} z
  = 0
  \label{prima}
  \\
  && \quad
  \hbox{\aat, for every $z\in V$,
    and for $\xi$ given by \Seconda}.
    \non
\Esist
\Accorpa\regu regua regub
\Accorpa\pbl regua prima

\Brem
\label{Sensosoluz}
We regard the initial/boundary value problem \pbl\ as an essentially integrodifferential Allen-Cahn equation
in the sole unknown~$\u$.
We note, in particular, that \eqref{prima} has a precise meaning, because 
$\xi^{1/2}\in\LQ2$ and $\u^{-1/2}\in\LQ\infty$
(at~least)
whenever $\u$ satisfies \eqref{regua}
and $\bar\emme\in\LQ2$.
\Erem

Here is our first result (the symbol $(\cpto)^-$ denotes the negative part).

\Bthm
\label{Wellposedness}
Assume that:
\Bsist
  &&\bar\emme \in \LQi\infty 
  \aand
  \emmem \in \LL1\Linfty;\quad \frac 1\uz + \frac 1{1-\uz} \in \Linfty, 
  \label{hpemme}
  \\
  && \uz \in \Hdue, \quad
  \dn\uz = 0 \quad \hbox{on~$\Gamma$},
  \aand
  \Delta\uz \in \Linfty .
  \label{highreg}
\Esist
Then, for every $T\in(0,+\infty)$,
problem \pbl\ has a unique solution.
Furthermore, 
\Beq
  \u \in \L p{\Wx{2,p}}
  \quad \hbox{for every $p<+\infty$}, \quad
  \dt\u \in \LQ\infty ,
  \aand
  \xi \in \LQ\infty .
  \label{piuregsoluz}
\Eeq
Finally, there exist constants $\umin,\umax\in(0,1)$
and $\ximax\geq0$ such~that
\Beq
  \umin \leq \u \leq \umax
  \aand
  \xi \leq \ximax \quad \aeQ;
  \label{stimesoluz}
\Eeq 
these constants can be chosen
independently of~$T$.
\Ethm

\Brem
\label{Strongsol}
Let $\Mz$ be an upper bound for $|\bar\emme|$. Assume that $\xi_*:=\inf\xiz$ is strictly positive, and take $M$ such that
$|\dt\u|\leq M$ \aeQ\
(see~\eqref{piuregsoluz}).
Then, relation \eqref{secondaloc}~ in the Appendix implies that $\xi(t)>0$ 
at least for $t<2\sqrt{\xi_*\umin}/(\kappa M^2+\Mz)$.
In such a case, the pair $(\u,\xi)$ 
yields ~a (local) solution $(\u,\mu,\eta)$
to the original problem in a strong sense.
\Erem

Our second result concerns the long-time \bhv\ of the solution~$\u$ to problem \pbl;
it ensures that the elements of the \omegalimit\ of every trajectory are steady states. 
To state this result properly, we have to describe the stationary problem
associated to~\pbl.
We introduce $\phii:\Omega\to[0,+\infty)$, by means of the following formula:
\Beq 
  \phii(x) := \lim_{t\to+\infty} \sqrt{\xi(x,t)}
  \quad \aaO ,
  \quad
  \hbox{where $\xi$ is given by \Seconda}
  \label{defphii}
\Eeq
(our next theorem 
ensures that such a limit actually exists, under form of a bounded function on~$\Omega$).
The stationary problem consists in finding
\Beq
  \ui \in V
  \aand
  \umin \leq \ui \leq \umax
  \quad \aeO
  \label{regui}
\Eeq
such that
\Beq
  \iO \nabla\ui \cdot \nabla z
  + \iO f'(\ui) \, z
  - \iO \frac \phii {\sqrt{\ui}} \, z
  = 0 
  \quad \hbox{for every $z\in V$}.
  \label{primai}
\Eeq
\Accorpa\pbli regui primai

\Bthm
\label{Longtime}
Under the assumptions of Theorem~\ref{Wellposedness},
let $\u$ be
the unique global solution to problem~\pbl.
Then, the limit~\eqref{defphii}
exists \aaO\ and $\phii\in\Linfty$.
Moreover,
the \omegalimit\ defined~by
\Beq
  \omega(\u)
  := \graffe{
    \uomega \in H : \ 
    \uomega = \limn \u(\tn) \ 
    \hbox{strongly in $H$ for some $\{\tn\}\nearrow+\infty$}
  }
  \label{defomega}
\Eeq
is non-empty, compact, and connected
in the strong topology of~$H$.
Finally, every element $\uomega\in\omega(\u)$
coincides with a solution $\ui$ 
to the stationary problem~\pbli.
\Ethm

\Brem
\label{Hpmaxmon}
One can wonder whether $\fu$ can be a more general potential from Convex Analysis.
Actually this is the case, 
since we may replace the monotone part $\fu'$
of $f'$ by a graph~$\alpha$.
Precisely, we may assume~that
\Bsist
  && \hbox{$\alpha$ is a maximal monotone graph in $\erre\times\erre$},
  \label{hpalphaa}
  \\
  &&  \hbox{with $D(\alpha)=(0,1)$,\quad $\lim_{r\searrow 0} \alpha^0(r) = - \infty,$ 
  \aand
  $\lim_{r\nearrow 1} \alpha^0(r) = + \infty$}
  \label{hpalphab},
\Esist
where $D(\alpha)$ stands for the domain of $\alpha$
and $\alpha^0(r)$ denotes the element of $\alpha(r)$
having minimum modulus for $r\in(0,1)$
(see, e.g.,~\cite[p.~28]{Brezis}).
Accordingly, $\fu$~is replaced by a convex l.s.c.\
function $\hat\alpha:\erre\to(-\infty,+\infty]$
such that $\partial\hat\alpha=\alpha$, so that $f:=\hat\alpha+\fd\geq0$.

When thinking of such a \generaliz ation,
we have to introduce a selection $\zeta$ of $\alpha(\u)$
with some regularity and we have to come up with a convenient replacement for the variational equation~\eqref{prima}.
Here is a suitably general formulation:
\Bsist
  && \textrm{find}\quad\zeta \in \LQ2
  \aand 
  \zeta \in \alpha(\u) \quad \aeQ,\quad\textrm{such that}
  \label{regzeta}
  \\
  && \kappa \iO \dt\u(t) \, v
  + \iO \nabla\u(t) \cdot \nabla v
  + \iO \bigl( \zeta(t) + \fd'(\u(t)) \bigr) \, v
  - \iO \bigl( \xi(t)/\u(t) \bigr)^{1/2} v
  = 0
  \qquad
  \non
  \\
  && \quad
  \hbox{\aat, for every $v\in V$,
    and for $\xi$ given by \eqref{seconda}}.
  \label{primagen}
\Esist
An analogous modification is due for the stationary problem~\pbli, that may be replaced
by the following problem:
\Bsist
  &&  \textrm{find}\quad\zetai \in \Ldue
  \aand 
  \zetai \in \alpha(\ui) \quad \aeO,\quad\textrm{such that}
  \label{regzetai}
  \\
  && \iO \nabla\ui \cdot \nabla z
  + \iO \bigl( \zetai + \fd'(\ui) \bigr) \, z
  - \iO \frac \phii {\sqrt{\ui}} \, z
  = 0 
  \quad \hbox{for every $z\in V$}.
  \label{primaigen}
\Esist
With such measures, Theorem~\ref{Wellposedness} can be extended
as far as existence and regularity are concerned.
Precisely, we can prove that there is a~global solution $(\u,\zeta)$, 
with $\zeta\in\LQi\infty$ and $\u$
satisfying the same regularity requirements and bounds
as in Theorem~\ref{Wellposedness}.
However, we cannot prove uniqueness.
Furthermore, Theorem~\ref{Longtime}
holds for every (possibly nonunique) global solution
to the \generaliz ed problem satisfying the same bounds as above. We sketch how to achieve such \generaliz ations
in the forthcoming Remarks~\futuro\ref{Maxmon} 
and~\futuro\ref{Maxmonbis}.
\Erem


\section{Existence and uniqueness}
\label{WELLPOSEDNESS}
\setcounter{equation}{0}

In this section, we prove Theorem~\ref{Wellposedness}. This is a rather complicated task, that we now delineate.
First of all, we show that every solution
satisfies the last part of the statement,
i.e., that some kind of maximum principle holds.
Then, we show that we can count a~priori on more regularity
than that specified in~\eqref{piuregsoluz}.
Finally, by looking for solutions satisfying 
such stronger properties, only, we prove 
existence and uniqueness using a fixed point argument.
In the preliminary steps, we find convenient to deal with an auxiliary problem.
In the whole section, it is understood that the assumptions
of Theorem~\ref{Wellposedness} are satisfied.

\step
Construction of the crucial constants

We find the constants 
$\umin$, $\umax$, and~$\ximax$,
noting that our procedure actually yields values
that do not depend on~$T$,
as stated in the last part of Theorem~\ref{Wellposedness}.
For convenience, we set
\Beq
  \Md := \sup_{r\in(0,1)} |\fd'(r)| ,
  \label{defMd}
\Eeq
and choose $\umin\in(0,1)$ and $\ximaxz\geq0$ so as to have
\Beq
  \fu'(\umin) \leq -\Md , \quad
  \uz \geq \umin \aand \xiz \leq \ximaxz \quad \aeO ,
  \label{defumin}
\Eeq
due to~\eqref{hpfb}, \eqref{hpdati}, and~\eqref{hpemme}.
Moreover, on accounting for the second of~\eqref{hpemme},
we define $\ximax\geq0$ as follows:
\Beq
  \sqrt\ximax := \sqrt\ximaxz + \frac 1{2\sqrt\umin} \, \norma{\emmem}_{\LL1\Linfty} \,.
  \label{defximax}
\Eeq
Finally, using the last of~\eqref{hpfb} and \eqref{hpdati} once more,
we choose $\umax\in(0,1)$ such~that
\Beq
  \fu'(\umax) - \frac{\sqrt\ximax}{\sqrt\umax} \geq \Md
  \aand
  \uz \leq \umax \quad \aeO .
  \label{defumax}
\Eeq

At this point, we are ready to prove 
the last part of Theorem~\ref{Wellposedness}.
At the same time, with a view towards the fixed 
point argument we are going to use
later on, we prepare some auxiliary material.
So, we show~that
\Beq
  \umin \leq \u \leq \umax
  \aand
  \xi \leq \ximax
  \quad \aeQ
  \label{principiomax}
\Eeq
for any solution $\u$ to the variational equation~\eqref{prima},
but with \Seconda\ replaced by something else.

\step
The auxiliary problem 

From the previous section, 
we see that~$\sqrt\xi$, rather than~$\xi$, 
plays the main role.
Hence, we define 
$\Phi:D(\Phi)\to\W{1,1}\Luno\cap\LQ\infty$ 
as follows.
We~set:
\Beq
 \DPhi := \Graffe{
  v \in \H1H, \quad
  v > 0 \quad \aeQ , \quad
  1/v \in \LQ\infty }
  \label{defDPhi}
\Eeq
and, for $v\in\DPhi$, we denote by $\Phi(v)$ 
the function $\phi$ given~by
\Beq
  \phi(x,t)
  = \sup \Graffe{w(x,t) :\ w\in\calF(v)}
  \quad \hbox{for $(x,t)\in Q_T$,}
  \label{genseconda}
\Eeq
where we have~set
\Bsist
  && \calF(v) 
  := \Bigl\{
  w\in\W{1,1}\Luno :\ w(0) = \sqxz, \quad w \geq 0 \quad \aeQ, 
  \non
  \\
  && \phantom{\calF(v) := \Bigl\{} \quad 
    \dt w = - \bigl( \kappa \, (\dt v)^2 + \bar\emme \bigr) / (2v^{1/2})
    \quad \hbox{a.e.\ where $w>0$}
  \Bigr\} .
  \label{defF}
\Esist
By arguing as in the Appendix, we see that
$\phi:=\Phi(v)$~is the square root of the maximal solution
to the Cauchy problem:
\Beq 
  \dt \xi 
  = - \Bigl(
       \bigl(
         \kappa \, (\dt v)^2 + \emme
       \bigr) / (2v^{1/2})
  \Bigr) \sqrt\xi
  \aand
  \xi(0) = \xiz , 
  \non
\Eeq
and it is \characteriz ed~by
\Beq
   \phi(x,t)
  = \sqrt{\xiz(x)} - \iot \astar(x,s) \ds ,
  \label{gensecondaa}
  \Eeq
where
\Beq
  \astar(x,s) := \frac{\kappa \, |\dt v(x,s)|^2 + \bar\emme(x,s)}{2 \, \sqrt{v(x,s)}}
  \quad \hbox{if} \quad \phi(x,t) > 0 , 
  \quad\ 
  \astar(x,s) := 0 \quad \hbox{otherwise}.
  \label{gensecondab}
\Eeq
Note that $\phi$ actually belongs
to $\W{1,1}\Luno\cap\LQ\infty$.
Then, the auxiliary problem is obtained as follows.
For a given $v\in\DPhi$,
we require that $\u$ satisfies:
\Bsist
  && \u \in \H1H \cap \C0V ,
  \quad \ 
  \u(0) = \uz , 
  \label{reguabis}
  \\
  && 0 < \u < 1 \quad \aeQ , \quad
  f'(u) \in \LQ2 ,
  \aand
  \u^{-1/2} \in \LQ2 ;
  \label{regubbis}
  \\
  && \kappa \iO \dt\u(t) \, z
  + \iO \nabla\u(t) \cdot \nabla z
  + \iO f'(\u(t)) \, z
  - \iO \frac {\phi(t)} {\sqrt{\u(t)}} \, z
  = 0
  \\
  && \quad
  \hbox{\aat\ and every $z\in V$,
    where $\phi=\Phi(v)$}.
     \non
  \label{primabis}
\Esist
\Accorpa\pblaux reguabis primabis
Therefore, problem \pbl\
is equivalent to the auxiliary problem, 
provided $v=\u$ and $\phi=\sqrt\xi$, and provided that 
some stronger regularity requirements are granted. 

\Blem
\label{Stimadalbasso}
Let $v\in\DPhi$ and assume that $\u$
satisfies \pblaux.
Then, we~have~that
\Bsist
  && \u \geq \umin \quad \aeQ .
  \label{princmin}
\Esist
In particular, this is true
if $v$ is any solution $\u$ to problem~\pbl.
\Elem

\Bdim
The proof we give is quite standard.
Let $g:\erre\to\erre$ be \Lip\ continuous, nondecreasing, and 
such that $g(r)<0$ for $r<\umin$ and $g(r)=0$ for $r\geq\umin$;
furthermore, let $G$ be the primitive of $g$ that vanishes at~$\umin$.
Now, we write~\eqref{primabis} at $t=s$ 
and test it by $z:=g(\u(s))$.
Then, we integrate over~$(0,t)$ with respect to~$s$,
where $t\in(0,T)$ is arbitrary.
We~have: 
\Beq
  \kappa \iO G(\u(t))
  + \intQt \nabla\u \cdot \nabla g(\u)
  + \intQt \tonde{f'(\u) - \phi/\squ} \, g(\u)
  = \kappa \iO G(\uz) .
  \label{perminmax}
\Eeq
The \rhs\ vanishes by~\eqref{defumin}
and the first two terms on the \lhs\ are nonnegative.
As to the third, we can replace $Q_t$
by its subset where $\u<\umin$.
Due to~\Hpf\ and~\accorpa{defMd}{defumin},
a.e.\ in such a subset we have that:
\Beq
  f'(\u) - \frac\phi\squ
  \leq f'(\u) 
  \leq \fu'(\umin) + \Md 
  \leq 0
  \aand
  g(\rho) \leq 0 ,
\Eeq
whence the nonnegativity of the third integral in \eqref{perminmax}. 
We conclude~that $G(\u(t))=0$ \aeO\ for every $t\in[0,T]$,
i.e., that $\u\geq\umin$ \aeQ.
\Edim

\Blem
\label{Stimaxi}
Assume $v\in \DPhi$ and $v\geq\umin$.
Then, we have that
\Beq
  \Phi(v) \leq \sqrt\ximax \quad \aeQ . 
  \label{stimaphi}
\Eeq
In particular, $\xi\leq\ximax$ \aeQ\
for every solution $(\u,\xi)$ to problem~\pbl.
\Elem

\Bdim
We set $\phi:=\Phi(v)$ for brevity
and stipulate that
\Beq
  \hbox{$\chi$ is the characteristic function
  of the subset of $Q_T$ where $\phi>0$}.
  \label{defchi}
\Eeq
Then, \accorpa{gensecondaa}{gensecondab}~yield
\Bsist
  && \phi(t)
  = \sqxz - \iot \chi(s) \, 
      \frac {\kappa \, |\dt v(s)|^2 + \emmep(s) - \emmem(s)}
      {2 v(s)}
  \ds
  \non
  \\
  && \leq \sqxz + \frac 12 \iot \chi(s) \, \frac{\emmem(s)}{\sqrt{v(s)}} \ds
  \leq \sqrt\ximaxz + \frac 1{2\sqrt\umin} \, \norma{\emmem}_{\LL1\Linfty}
  = \sqrt\ximax \,.
  \non
\Esist
In particular, $\xi\leq\ximax$ if $v=\u$
and $(\u,\xi)$ solves problem \pbl,
because the inequality $\u\geq\umin$ 
has been already proved
for every solution.
\Edim

\Blem
\label{Stimadallalto}
Assume $v\in\DPhi$ and $v\geq\umin$ \aeQ\
and let $\u$ satisfy \pblaux.
Then, we have~that
\Beq
  \u \leq \umax \quad \aeQ . 
  \label{princmax}
\Eeq
In particular, this is true
if $v$ is any solution $\u$ to problem~\pbl.
\Elem

\Bdim
Let $g:\erre\to\erre$ be \Lip\ continuous, nondecreasing,
such that $g(r)>0$ for $r>\umax$ and $g(r)=0$ for $r\leq\umax$;
and let $G$ be the primitive of $g$ that vanishes at~$\umax$.
Then, \eqref{perminmax}~holds with the new $g$ and~$G$;
once again, the only trouble comes from the third term on the \lhs.
However, (a.e.)~in the subset of $Q_t$ where $\u>\umax$, we have~that
\Beq
  g(\u) \geq 0
  \aand
  f'(\u) - \frac \phi\squ
  \geq \fu'(\u) - \frac {\sqrt\ximax} \squ + \fd'(\u)
  \geq \fu'(\umax) - \frac {\sqrt\ximax} {\sqrt\umax} - \Md
  \geq 0,
\Eeq
thanks to the estimate~\eqref{stimaphi}
and~the definition of $\umax$ given by~\eqref{defumax}.
We conclude~that $G(\u(t))=0$ \aeO\ for every $t\in[0,T]$,
i.e., that $\u\leq\umax$ \aeQ.
\Edim

\Brem
\label{Ultimafrase}
We note that the last sentence of Theorem~\ref{Wellposedness}
is completely proved.
\Erem

\Blem
\label{HigherLpreg}
Let $v\in\DPhi$ be such that $\umin\leq v\leq\umax$ \aeQ,
and let $\u$ satisfy \pblaux.
Then, we have that
\Bsist
  && \norma{\dt\u}_{\LQ p}
  \leq \Rp
  \quad \hbox{for every $p\in(1,+\infty)$},
  \label{stimaLp}
  \\
  && \norma\u_{\L p{\Wx{2,p}}}
  \leq \Rp'
  \quad \hbox{for every $p\in(1,+\infty)$},
  \label{stimaWp}
\Esist
where the constants $\Rp$ and $\Rp'$ 
depend only on the structure
of our problem, the initial data $\uz$ and~$\xiz$, and~$p$.
In particular, all this is true
if $v$ is any solution $\u$ to problem~\pbl.
\Elem

\Bdim
Observe that
\Beq
  \iO \dt\u(t) z + \iO \nabla\u(t) \cdot \nabla z 
  = \iO F(t) z
  \aand
  \u(0) = \uz
  \non
\Eeq  
\aat, for every $z\in V$,
and for $F:=\phi/\squ-f'(\u)$.
Owing to~\eqref{princmin}, \eqref{stimaphi}, \eqref{princmax},
and to our assumptions on~$f$ and~$\uz$, we~have~that
\Beq
  |F| \leq \frac{\sqrt\ximax}{\sqrt\umin}
  + \sup_{\umin\leq r\leq\umax} |f'(r)|
  \quad \aeQ
  \aand
  \norma\uz_{\Linfty} \leq \umax.
  \label{perdefRp}
\Eeq
Therefore, we can apply the general
$L^p$-regularity theory
(see, e.g., \cite[Thm.~9.1, p.~341]{LSU})
and deduce~that
\accorpa{stimaLp}{stimaWp} hold.
Moreover, thanks to the above lemmas,
the first of~\eqref{piuregsoluz} is proved.
\Edim

\Blem
\label{StimaLinfty}
Let $v\in\DPhi$ be such that $\umin\leq v\leq\umax$ \aeQ,
and let $\u$ satisfy \pblaux.
Moreover, assume~that
\Beq
  \norma{\dt v}_{\LQ p} \leq \Rp
  \quad \hbox{for every $p\in(1,+\infty)$},
  \label{hpdtLp}
\Eeq
with the same $\Rp$ as in~\eqref{stimaLp}.
Then, 
\Beq
  \norma{\dt\u}_{\LQ\infty}
  \leq \Rinfty ,
  \label{stimaLinfty}
\Eeq
for some constant $\Rinfty$ depending only on the structure
of our problem and the initial data $\uz$ and~$\xiz$.
In particular, this is true if $v$ is any solution $\u$ to problem~\pbl.
\Elem

\Bdim
Here we use the stronger regularity assumption~\eqref{highreg}.
We proceed formally, because the argument that would make the calculation rigorous
is quite standard. We differentiate~\eqref{primabis} 
with respect to time and~obtain:
\Beq
\kappa \iO \dt u(t) \, z + \iO \nabla u(t) \cdot \nabla z + \iO u(t) \, z
  = \iO F(t) \, z
  \ \
 \hbox{\aat\ and every $z\in V$,}
 \label{derivata}
\Eeq
where we have set set $u:=\dt\u$ and 
\Beq
  F := \dt\u - f''(\u) \dt\u 
  - \frac 12 \, \phi \, \u^{-3/2} \dt\u
  - \frac 12 \, \chi \, \bigl( \kappa |\dt v|^2 + \bar\emme \bigr) (v\u)^{-1/2} ,
  \non
\Eeq
with $ \chi$ as in~\eqref{defchi}.
Now, the initial value of $u$ is known explicitely
through the differential equation for~$\u$, and
is given by the formula:
\Beq
  \kappa u(0) = \Delta\uz - f'(\uz) + \sqrt\xiz/\sqrt\uz \,.
  \non
\Eeq
Hence, $u(0)$ belongs to $\Linfty$ and its $L^\infty$-norm
is estimated by a known constant, due to~\eqref{highreg}.
Now, we observe that the following estimate holds 
(see \cite[Thm.~7.1, p.~181]{LSU}):
\Beq
  \norma u_{\LQ\infty}
  \leq C_q \max \graffe{ \norma{u(0)}_{\Linfty} + \norma F_{L^q(Q_T)}},
  \label{stimaKS}
\Eeq
where $q$ and $C_q$ are as follows:
$q$~is any real number satisfying $q>r_N^2/(r_N-1)$,
where $r_N=\max\graffe{2,1+(N/2)}$;
$C_q$ is a constant depending only on $\Omega$ and~$q$.
Hence, we are led to find an $L^q$-estimate of $F$ for some $q$
satisfying the above inequality.
Actually, we can choose anyone of such values of~$q$
and get an estimate using
the previous lemma and assumptions~\eqref{hpemme} and~\eqref{hpdtLp}.
If we denote by $C$ a suitable constant
that could be computed in terms of
our assumptions on the problem structure and the constants
$\umin$, $\umax$, and~$\ximax$, we have indeed that
\Bsist
  && \norma F_{\LQ q}
  \leq C \tonde{
    \norma{\dt\u}_{\LQ q}
    + \norma{(\dt v)^2}_{\LQ q}
    + \norma{\bar\emme}_{\LQ q}
  }
  \non
  \\
  && \leq C ( R_q + R_{2q}^2 + \norma{\bar\emme}_{\LQ q}),
  \non
\Esist
for any $q\in(1,+\infty)$. 
Hence, \eqref{stimaKS}~provides a value of $\Rinfty$
satisfying~\eqref{stimaLinfty}.
\Edim

\Brem
\label{Ogniregolarita}
We note, in particular,
that now all the conditions listed in~\eqref{piuregsoluz} are proved.
\Erem

\step
Towards the fixed point argument

Due to the above lemmas,
we can confine ourselves to look for solutions
satisfying all the bounds we have proved,
i.e., belonging to the set $\calX$ defined here below:
\Beq
  \calX :=
  \graffe{ v\in\calV: \ \umin \leq v \leq \umax \ \aeQ , \
    \norma{\dt v}_{\LQ p} \leq \Rp \ \,
    \forall \, p\in(1,+\infty]
   },
   \qquad
  \label{defX}
\Eeq
where the constants involved are of course the same as before,
and where
\Beq
  \calV :=  \H1H \cap \C0V .
  \label{defV}
\Eeq
We regard $\calX$ as a metric subspace of the Banach space~$\calV$, i.e.,
we consider the metric $d$ on $\calX$ defined by
\Beq
  d^2(u,v) 
  := \intQ |\dt u - \dt v|^2 + \sup_{t\in[0,T]} \normaV{u(t)-v(t)}^2
  \quad \hbox{for $u,v\in\calX$}.
  \label{defd}
\Eeq
Thus, $\calX$ is complete: each one of the conditions specified 
in the definition~\eqref{defX} of $\calX$ 
defines a closed subset of~$\calV$,
and such is $\calX$, being the intersection
of a family of closed sets. 
Next, for $v\in\calX$, we observe that $v\in\DPhi$,
and we define the mapping $\Psi:\calX\to\calX$ as follows:
\Beq
  \hbox{$\Psi(v)$~is the unique solution~$\u$
  to problem \pblaux}.
  \label{defPsi}
\Eeq
It is clear that, for every $v\in\calX$,
the auxiliary problem \pblaux\
has a unique solution~$\rho$.
Indeed, if we set $\phi:=\Phi(v)$,
we see that, \aaQ, the nonlinear function
$r\mapsto f'(r)-\phi(x,t)/\sqrt r$,
defined for
$\umin\leq r\leq\umax$,
is \Lip\ continuous;
moreover, $\phi$~is bounded.
Hence, the standard theory for regular
parabolic equations yields
existence and uniqueness
and we conclude that $\Psi(v)$ is well defined.
Furthermore,
the precise choice of all constants in~\eqref{defX}
exactly ensures that $\Psi(v)$ belongs to~$\calX$,
due to the previous lemmas,
whence $\Psi$ actually maps $\calX$ into itself.
Finally, it is clear that an element of $\calX$
is a solution to \pbl\ if and only if
it is a fixed point for~$\Psi$.
Hence, it suffices to prove that $\Psi^k:=\Psi\circ\dots\circ\Psi$
($k$~times) is~a contraction for $k$ large enough.
The rest of the proof is devoted to do~that.
A key point of our argument is the following lemma. 

\Blem
\label{ContrazODE}
For $i=1,2$, pick $a_i\in\Lloc$,
and let $y_i$ be the maximal solution
to the Cauchy problem~\eqref{ode} with $a=a_i$.
Then, we~have
\Beq
  \sup_{s\in[0,t]} \bigl| \sqrt{y_1(s)} - \sqrt{y_2(s)} \bigr|
  \leq \iot |a_1(s) - a_2(s)| \dtau
  \quad \hbox{for every $t\geq0$}.
  \label{contrazODE}
\Eeq
\Elem

\Bdim
We first prove~that
\Beq
  \frac d{dt} \, \bigl| \sqrt{y_1} - \sqrt{y_2} \bigr|
  \leq |a_1 - a_2 |
  \quad \hbox{a.e.\ in $(0,+\infty)$}.
  \label{pratelli}
\Eeq
To this end,
we choose everywhere defined representatives of~$a_i$,
set $\phi_i:=\sqrt{y_i}$ for brevity,
and notice that the functions
$\phi_i$ and $|\phi_1-\phi_2|$
are locally absolutely continuous on~$[0,+\infty)$.
Therefore, there exists an exceptional set 
$E\subset[0,+\infty)$ having zero Lebesgue measure,
such that the derivatives at~$t$ of the above 
functions exist for every $t\in(0,+\infty)\meno E$.
Moreover, we have that $\phi'_i(t)=-a_i(t)$ if $\phi_i(t)>0$.
Let us fix a point~$t$ outside of~$E$
and prove~\eqref{pratelli} at~$t$.
We distinguish three cases.
In the first one, we have
$\phi_i(t)>0$ for $i=1,2$.
Then, $|\phi_1-\phi_2|=\pm(\phi_1-\phi_2)$
in a \nbh\ of $t$ and
$\phi_i'(t)=-a_i(t)$ for $i=1,2$,
whence
\Beq
  |\phi_1-\phi_2|'(t)
  = \pm(\phi_1-\phi_2)'(t)
  = \mp(a_1-a_2)(t)
  \leq |a_1(t) - a_2(t)|.
  \non
\Eeq
Assume now $\phi_i(t)=0$ for $i=1,2$.
Then, $\phi'_i(t)=0$ (see~\eqref{soluzodeb})
and~the desired inequality trivially follows.
In the last case, we have, e.g.,
$\phi_1(t)>0=\phi_2(t)$.
We derive that $\phi'_1(t)=-a_1(t)$.
On the other hand, as $\phi_2(t)=0$,
we see that \eqref{perpratelli} implies
$a_2(t)\geq0$.
Therefore, noting that 
$\phi_1-\phi_2>0$ in a \nbh\ of~$t$,
we deduce~that
\Beq
  |\phi_1 - \phi_2|'(t)
  = (\phi_1 - \phi_2)'(t)
  = \phi'_1(t) - \phi'_2(t)
  \leq - a_1(t) + a_2(t)
  \leq |a_1(t) - a_2(t)|,
  \non
\Eeq
and \eqref{pratelli} is completely proved.
Now, we derive~\eqref{contrazODE}.
We fix $t\geq0$ and $s\in[0,t]$.
Then, \eqref{pratelli}~yields
\Beq
  \bigl| \sqrt{y_1(s)} - \sqrt{y_2(s)} \bigr|
  = \ios \bigl| \sqrt{y_1} - \sqrt{y_2} \bigr|'(\tau) \dtau
  \leq \ios |a_1(\tau) - a_2(\tau) | \dtau
  \leq \iot |a_1(\tau) - a_2(\tau) | \dtau ,
  \non
\Eeq
and \eqref{contrazODE} immediately follows.
\Edim

\step
Conclusion 

We take $v_1,v_2\in\calX$ and set for convenience
$\phi_i:=\Phi(v_i)$ and $\u_i:=\Psi(v_i)$ for $i=1,2$.
Then, we write the equality in \eqref{primabis}
for $v=v_i$ and test the difference by $\dt(\u_1-\u_2)$.
Then, we integrate over~$(0,t)$ for an arbitrary $t\in(0,T)$
and add the same integral to both sides for convenience.
If we set $\u:=\u_1-\u_2$ and use a similar notation
for $\phi$ and~$v$, we obtain:
\Beq
  \kappa \intQt  |\dt\u|^2
  + \frac 12 \normaV{\u(t)}^2
  = \intQt \tonde{\u - f'(\u_1) + f'(\u_2)} \dt\u
  + \intQt \tonde{\phi_1 \u_1^{-1/2} - \phi_2 \u_2^{-1/2}}
  \dt\u 
  \label{percontraz}
\Eeq
We now pass to estimate the \rhs\ of the last relation. 
In order to simplify the notation,
we use the same symbol $c$
for different constants (even~in the same formula)
that only depend on the problem structure, the data, 
and~$T$; $\cd$ denotes  those such constants that 
depend, in addition,  on the parameter $\delta\in(0,1)$.
As the functions $f$ and $(\cpto)^{\pm1/2}$ 
are \Lip\ continuous on the interval~$[\umin,\umax]$, 
and as estimate~\eqref{stimaphi}
holds for both~$\phi_i$,
we have that
\Bsist
  && \intQt \tonde{\u - f'(\u_1) + f'(\u_2)} \dt\u
  + \intQt \tonde{\phi_1 \u_1^{-1/2} - \phi_2 \u_2^{-1/2}}
  \dt\u 
  \non
  \\
  && \leq c \intQt |\u| \, |\dt\u|
  + \intQt \phi_1 \, |\u| \, |\dt\u|
  + \intQt \u_2^{-1/2} |\phi| \, |\dt\u|
  \non
  \\
  && \leq \delta \intQt |\dt\u|^2 
  + \cd \intQt |\u|^2
  + \cd \intQt |\phi|^2 ,
  \non
\Esist
for every $\delta\in(0,1)$. By combining with \eqref{percontraz} and choosing 
$\delta$ small enough, we easily deduce that
\Beq
  \intQt  |\dt\u|^2
  + \normaV{\u(t)}^2
  \leq c \iot \normaH{\u(s)}^2 \ds
  + c \iot \normaH{\phi(s)}^2 \ds.
  \non
\Eeq
On the other hand, we can compare \Seconda,
\accorpa{soluzode}{soluzodeb}, 
and \accorpa{defDPhi}{gensecondab}
and apply Lem\-ma~\ref{ContrazODE}.
Thus, for every $s\in(0,t)$ and \aeO, we have~that
\Bsist
  && |\phi(s)| = |\phi_1(s) - \phi_2(s)|
  \non
  \\
  && \leq \frac 12 
  \ios \left|
    v_1^{-1/2}(\tau) \bigl( \kappa \, |\dt v_1(\tau)|^2 + \bar\emme(\tau) \bigr)
    - v_2^{-1/2}(\tau) \bigl( \kappa \, |\dt v_2(\tau)|^2 + \bar\emme(\tau) \bigr)
  \right| \dtau
  \non
  \\
  && \leq \frac 12 \ios |v_1^{-1/2}(\tau) - v_2^{-1/2}(\tau)| \, 
    \bigl( \kappa \, |\dt v_1(\tau)|^2 + |\bar\emme(\tau)| \bigr) \dtau
  \non
  \\
  && \quad {}
  + \frac \kappa 2 \ios v_2^{-1/2}(\tau)
    \bigl|
      |\dt v_1(\tau)|^2 - |\dt v_2(\tau)|^2
    \bigr|
  \non
  \\
  && \leq c \ios \bigl( 1 + |\bar\emme(\tau)| \bigr) |v(\tau)| \dtau
  + c \ios |\dt v(\tau)| \dtau
  \leq c \ios |v(\tau)| \dtau
  + c \ios |\dt v(\tau)| \dtau ,
  \non
\Esist
since $\umin\leq v_i\leq\umax$.
Hence,
\Bsist
  && \normaH{\phi(s)}^2
  \leq c \norma v_{\LQs2}^2
  + c \norma{\dt v}_{\LQs2}^2 
  \leq  c \intQs |\dt v|^2 
  + c \norma v_{\Cs0V}^2 .
  \non
\Esist
Therefore, we deduce that
\Beq
  \intQt |\dt\u|^2
  + \normaV{\u(t)}^2
  \leq c \iot \normaV{\u(s)}^2 \ds
  + c \iot 
  \tonde{
    \int_{Q_s} |\dt v|^2
    + \norma v_{\Cs0V}^2 
  } \ds ,
  \non
\Eeq
and the Gronwall lemma easily yields:
\Beq
  \intQt |\dt\u|^2
  + \norma\u_{\Ct0V}^2
  \leq \Cstar \iot 
  \tonde{
    \int_{Q_s} |\dt v|^2
    + \norma v_{\Cs0V}^2 
  } \ds
  \quad \hbox{for every $t\in[0,T]$,}
  \non
\Eeq
with a precise constant $\Cstar$.
This means~that
\Bsist
  && \intQt |\dt \bigl( \Psi^k(v_1) - \Psi^k(v_2) \bigr)|^2
  + \norma{\Psi^k(v_1) - \Psi^k(v_2)}_{\Ct0V}^2
  \non
  \\
  && \leq \frac{\Cstar^k}{(k-1)!} \iot s^{k-1}
  \tonde{
    \int_{Q_s} |\dt (v_1 - v_2)|^2
    + \norma{v_1-v_2}_{\Cs0V}^2 
  } \ds ,
  \non
\Esist
for every $t\in[0,T]$ with $k=1$.
Arguing by induction on~$k$,
it is \sfw\ to prove that the above inequality
holds for every~$k\geq1$.
We conclude~that
\Beq
  d^2(\Psi^k(v_1) - \Psi^k(v_2))
  \leq \frac{\Cstar^k \, T^k}{k!} \, d^2(v_1,v_2)
  \quad \hbox{for every integer $k\geq1$,}
  \label{contraz}
\Eeq
whence, $\Psi^k$ is a contraction on $\calX$ 
for $k$ large enough.
This completes the proof.

\Brem
\label{Maxmon}
Here, we breafly sketch how to modify
the above proof in order to achieve
the \generaliz ation mentioned in 
Remark~\ref{Hpmaxmon},
where \eqref{prima} was replaced 
by~\accorpa{regzeta}{primagen}.
Accordingly, we consider the auxiliary problem
obtained by assuming~\eqref{regzeta} for $\zeta$
and modifying \eqref{primabis} as~follows:
\Bsist
  && \kappa \iO \dt\u(t) \, z
  + \iO \nabla\u(t) \cdot \nabla z
  + \iO (\zeta + \fd'(\u(t))) \, z
  - \iO \bigl( \xi(t)/\u(t) \bigr)^{1/2} z
  = 0
  \\
  && \quad
  \hbox{\aat\ and every $z\in V$,
    where $\xi$ is given by \eqref{genseconda}}.
     \non
  \label{primagenbis}
\Esist
However, it is convenient to consider
approximating problems as~well.
Precisely, we write $\aeps(\u)$
instead of $\zeta$ in equations~\eqref{primagen}
and~\eqref{primagenbis},
where $\aeps$ denotes
the Yosida \regulariz ation of~$\alpha$
at level~$\eps\in(0,1)$
(see, e.g.~\cite[p.~28]{Brezis}).
We will speak of
the $\eps$-\generaliz ed problem and
of the $\eps$-auxiliary problem, respectively.
Now, we briefly show how to obtain first
uniform bounds as in the previous lemmas,
and then existence, for the \generaliz ed problem.
As far as the choice of the crucial constants
is concerned, in conditions~\eqref{defumin} 
and~\eqref{defumax}
we replace $\fu'$ and $\Md$
by $\alpha^0$ and~$\Md+1$, respectively.

As $|\aeps(r)|\leq|\alpha^0(r)|$,
and as $\aeps(r)$ converges to~$\alpha^0(r)$
when $\eps$ tends to zero for every $r\in(0,1)$,
we see that the inequalities
\Beq
  \aeps(\umin) \leq -\Md 
  \aand
  \aeps(\umax) - \frac {\sqrt{\ximax}} {\sqrt{\umax}}
  \geq \Md 
  \non
\Eeq
hold true for $\eps$ small enough.
Therefore, it is easy to see that
Lemmas~\ref{Stimadalbasso}--\ref{Stimadallalto}
still hold for each $\eps$-problem,
i.e., that the a~priori bounds~\eqref{principiomax}
are fulfilled in the new situation.

As to Lemma~\ref{HigherLpreg}, we observe~that
\Beq
  \alpha^0(r') \leq \aeps(r) \leq \alpha^0(r'')
  \quad \hbox{whenever $0<r'<\umin\leq r\leq\umax<r''<1$} .
  \label{perstimazeta}
\Eeq
Therefore, bounds~\eqref{principiomax}
imply uniform bounds for $\zeta=\aeps(\u)$
that can play the role of~\eqref{perdefRp}.
Thus, we see that \eqref{stimaLp} holds
for the solution to each of the $\eps$-problems
uniformly with respect~to~$\eps$.

As far as the \generaliz ation 
of Lemma~\ref{StimaLinfty} is concerned,
time differentiation is allowed
for $\eps$-problems.
However, we replace~\eqref{derivata}
by the following equality
\Beq
  \kappa \iO \dt u(t) \, z
  + \iO \nabla u(t) \cdot \nabla z
  + \iO u(t) \, z
  + \iO \aeps'(\u) \, u(t) \, z
  = \iO F(t) \, z, 
  \non
\Eeq
and we modify the previous definition of $F$
by writing $\fd''$ in place of~$f''$.
In other words, we replace
the term $-\fu''(\u)\,\dt\u$ on the \rhs\
by the term $\aeps'(\u)\,u$ on the \lhs.
As a possible technique for~\eqref{stimaKS}
consists in testing the above equation
by $z:=Z(u)$, where $Z$ is monotone and vanishes at~$0$
in~order to get recursive $L^p$-estimates
by a Moser's argument,
and as the integral involving $\aeps'$ 
gives a nonnegative contribution in such a procedure,
the type-\eqref{stimaLinfty} bound one finds 
for the solutions to our $\eps$-problems
is uniform with respect to~$\eps$.
All this ensures that the definition of~$\calX$
can be done in the present case, and it is independent of~$\eps$.
However, the constant $\Cstar$ we find
in applying the fixed point argument
does depend on~$\eps$.
Nevertheless, this is enough
to conclude for the existence of a unique solution
to the approximating problem
among the functions $\u$ belonging to~$\calX$.
Since uniqueness 
among all solutions can be proved
for every fixed $\eps>0$,
we see that our argument constructs 
a global solution $\ueps$
defined in the whole of~$[0,+\infty)$.
Moreover, for every fixed~$T$, 
$\ueps$~satisfies a number of
a~priori estimates 
uniformly with respect to~$\eps$,
and \eqref{stimesoluz}
and \eqref{perstimazeta} imply that
a uniform $L^\infty$-estimate holds for $\zetaeps:=\aeps(\ueps)$.
Therefore, modulo standard arguments,
we see that $(\ueps,\zetaeps)$ sequentially converges
in the proper topology to a pair $(\u,\zeta)$
and that $(\u,\zeta)$ is
a global solution of the \generaliz ed problem.
Unfortunately, in the new situation the previous uniqueness proof
does not work, because \eqref{primagen} cannot be differentiated 
with respect to time.
\Erem


\section{Long-time \bhv}
\label{LONGTIME}
\setcounter{equation}{0}

In this section, we prove Theorem~\ref{Longtime}.
Hence, we choose data satisfying the prescribed conditions
and pick the corresponding unique global solution~$(\rho,\xi)$.
Our proof is organized as follows.
The first lemma establishes the first assertion of the theorem,
i.e., that the function $\phii$
given by~\eqref{defphii} is well-defined and bounded.
In the next one, we derive some new a~priori estimates,
which ensure that the \omegalimit\
$\omega(\u)$ given by~\eqref{defomega}
has the desired properties.
Finally, we conclude by the announced 
\characteriz ation of $\omega(\u)$.
We find it convenient to set:
\Beq
  \phi := \sqrt\xi 
  \aand
  \hbox{$\chi :={}$the characteristic function
  of the subset of $Q_\infty $ where $\phi >0$,}
  \label{defphichi}
\Eeq
and we notice that the function $\phii$ to be studied
is the pointwise limit of~$\phi$ as time goes to infinity.
Furthermore, we recall that estimates 
\eqref{princmin} and \eqref{princmax} hold for~$\u$.

\Blem
\label{Defphii}
The limit~\eqref{defphii}
is well defined and $\phii$ is bounded.
\Elem

\Bdim
Equations \accorpa{gensecondaa}{gensecondab}
and definitions \eqref{defphichi}~yield:
\Bsist
  && \phi(t)
  = \sqxz - \iot \chi(s) \, 
      \frac {\kappa \, |\dt\u(s)|^2 + \emmep(s) - \emmem(s)}
      {2 \u(s)}
  \ds
  = \lambda_-(t) - \lambda_+(t)
\Esist
where we have set
\Beq
   \lambda_-(x,t) := \sqrt{\xiz(x)} + \iot \chi(x,s) \, \frac {\emmem(x,s)} {2\sqrt{\u(x,s)}} \ds
  \non
\Eeq
and
\Beq 
\lambda_+(x,t) := \iot 
      \chi(x,s) \, \frac {\kappa \, |\dt\u(x,s)|^2 + \emmep(x,s)} {2\sqrt{\u(x,s)}}
      \ds
  \non
\Eeq
\aaQi.
To prove the assertion, it suffices to show that $\lambda_\pm$ are bounded 
and that $\lambda_\pm(x,t)$ are convergent 
as $t$ tends to infinity \aaO.
We recall that $\phi$ is nonnegative.
Hence, in view also of the last of~\eqref{hpemme}, we have that
\Beq
  0 \leq \lambda_+(x,t)
  \leq \lambda_-(x,t)
  \leq \norma\xiz_{\Linfty}^{1/2}
  + \frac{\norma\emmem_{\L1\Linfty}}{2\sqrt{\umin}},
  \quad \aaQi ,
  \label{lambdabdd}
\Eeq
so that $\lambda_\pm$ are bounded.
On the other hand, it is clear that
both $\lambda_-$ and $\lambda_+$
are non-decreasing with respect to time, so that their
convergence is ensured.
\Edim

\Blem
\label{Nuovestime}
We~have that:
\Bsist
  && \phi \leq \sqrt{\ximax}
  \quad \aeQi 
  \aand
  \intQi |\dt\phi| < + \infty \, ;
  \label{stimephi}
  \\
  && \u \in \LL\infty V 
  \aand
  \intQi |\dt\u |^2 < +\infty \, .
  \label{stimerho}
\Esist
\Elem

\Bdim
The first of~\eqref{stimephi}
is a consequence of Lemma~\ref{Stimaxi}.
To prove the second one,
we observe that the calculation in the previous proof~yields:
\Beq
  \iO \lambda_+(t)
  \leq \iO \lambda_-(t)
  \leq c \, |\Omega|
  \quad \hbox{for every $t>0$,}
\Eeq
where $c$ is the \rhs\ of~\eqref{lambdabdd}
and $|\Omega|$ is the measure of~$\Omega$.
This clearly implies~that
\Beq
  \intQi \chi \, |\dt\u |^2 < +\infty 
  \aand
  \intQi \chi \, \bar\emme^\pm < +\infty,
  \label{chidt}
\Eeq
because $\u$ is bounded from below.
On the other hand, we have that
\Beq
  \dt\phi = - \chi \, \frac {\kappa |\dt\u|^2 + \bar\emme} {2 \sqrt\u},
  \label{dtphi}
\Eeq
whence immediately 
\Beq
  \intQt |\dt\phi|
  \leq \intQt \chi \, \frac {\kappa |\dt\u|^2 + \emmep + \emmem} {2\sqrt\u}
  = \iO \Bigl( \lambda_+(t) + \lambda_-(t) - \sqxz \Bigr)
  \leq 2 c \, |\Omega|
  \non
\Eeq
for every~$t>0$ (with the same meaning of $c$ as before);
the second of~\eqref{stimephi} follows.
To prove~\eqref{stimerho},
we formally test~\eqref{prima} by~$\dt\u$ 
and integrate over $(0,t)$.
Then, we perform an integration by parts
with respect to time, and~obtain:
\Bsist
  && \kappa \intQt |\dt\u|^2
  + \intQt 2\sqrt\u \, \dt\phi
  + \iO \tonde{
    \frac 12 \, |\nabla\u(t)|^2 + f(\u(t))
  } + 2 \iO \sqrt{\uz} \, \sqrt{\xiz}
  \non
  \\
  && = \iO \tonde{
    \frac 12 \, |\nabla\uz|^2 + f(\uz) 
  } 
  + 2 \iO \sqrt{\u(t)} \, \phi(t)
  \non
\Esist
for every $t>0$.
By accounting for~\eqref{dtphi},
we see that the above equality can be written as follows:
\Bsist
  && \kappa \intQt (1 - \chi) \,|\dt\u|^2
  + \iO \tonde{
    \frac 12 |\nabla\u(t)|^2 + f(\u(t))
  } + 2 \iO \sqrt{\uz} \, \sqrt{\xiz}
  \non
  \\
  && = \iO \tonde{
    \frac 12 \, |\nabla\uz|^2 + f(\uz) 
  } 
  + 2 \iO \sqrt{\u(t)} \, \phi(t)
  + \intQt \chi \, \bar\emme .
  \label{pernuova}
\Esist
Thus, all the terms on the \lhs\ 
of~\eqref{pernuova} are nonnegative.
Moreover, the \rhs\ is bounded, since
$\u$ and $\phi$ are bounded 
and $\chi\emme\in\LQi1$ by
the second of~\eqref{chidt}.
We immediately deduce that
\Beq
  \intQi (1-\chi) \, |\dt\u|^2 < +\infty,
  \non
\Eeq
i.e., that the first of~\eqref{stimerho} holds.
Moreover, on recalling the first of~\eqref{chidt},
we see that the second of~\eqref{stimerho}
holds as~well.
\Edim

\step
Conclusion of the proof

Properties~\eqref{stimerho} imply, in particular,
that $\u$ is a bounded weakly continuous
$V$-valued function on~$[0,+\infty),$
due to the compact embedding $V\subset H$.
Therefore, the set $\omega(\u)$ is a non-empty
compact subset of $H$.  Actually, $\omega(\u)$ is also connected, due
to the continuity of $\u$ from $[0,+\infty)$ to $H$ and
to a standard argument from the theory of dynamical systems (see,
for instance, \cite[p.~12]{Haraux}).
Then, the first properties of the \omegalimit\
stated in the theorem follow.
It remains for us to \characteriz e~$\omega(\u)$.

Let $\uomega\in\omega(\u)$ and let $\graffe{\tn}$
be a diverging sequence such that
$\graffe{\u(\tn)}$ converges to $\uomega$
strongly in~$H$.
Then, we define $\un$ and $\phin$ by the following formulas:
\Beq
  \un(t) := \u(t+\tn)
  \aand
  \phin(t) := \phi(t+\tn)
  \quad \hbox{for $t\geq0$,}
  \non
\Eeq
and we consider their weak limits
on a fixed bounded interval~$(0,T)$
(e.g.,~for $T=1$).
Our aim is to prove that
such limits do not depend on time
and furnish a steady state.
First of all, we notice~that
\Beq
  \kappa \intQ \dt\un \, z
  + \intQ \nabla\un \cdot \nabla z
  + \intQ f'(\un) \, z
  - \intQ \phin \, \un^{-1/2} \, z
  = 0
  \quad \hbox{for every $z\in V$,}
  \label{priman}
\Eeq
as one immediately sees by integrating~\eqref{prima}
over~$(\tn,\tn+T)$.
Next, we observe~that
\Beq
  \norma\un_{\L\infty V} \leq \norma\u_{\LL\infty V} 
  \aand
  \norma\phin_{\LQ\infty} \leq \norma\phi_{\LQi\infty}
  \label{perlongtime}
\Eeq
for every $n$.
Moreover, as $\phi(t)$ converges pointwise
to $\phii$ as $t$ tends to infinity,
we easily see that the whole sequence $\graffe{\phin}$
converges to $\phii$ weakly-star in~$\L\infty H$.
Furthermore, 
$\dt\phi$ belongs to~$\LQi1$ by~\eqref{stimaphi}.
By accounting for the second of~\eqref{perlongtime} as well,
we deduce~that
\Beq
  \phin \to \phii
  \quad \hbox{strongly in $\C0{\Lx p}$ for any $p<+\infty$}.
  \label{limphin}
\Eeq
Considering now \eqref{perlongtime} once more,
we infer~that, to within subsequences,
\Beq
  \un \to \ui
  \quad \hbox{weakly in $\L\infty V$}
  \aand
  \dt\un \to 0
  \quad \hbox{strongly in $\L2H$}
  \label{limun}
\Eeq
for some $\ui$. We deduce~that 
(see~\cite[Sect.~8, Cor.~4]{Simon})
\Beq
  \un \to \ui 
  \quad \hbox{strongly in~$\C0H$,}
  \label{convforte}
\Eeq
whence, in particular, $\un(0)\to\u(0)$ strongly in~$H$.
Moreover, \eqref{limun}~also imply that
$\ui$ must be a constant,
namely, $\ui(t)=\ui(0)$ for every $t\in[0,T]$.
On the other hand, $\un(0)\to\uomega$ by assumption.
We conclude that $\ui(t)=\uomega$ for every $t\in[0,T]$.
Furthermore, as $f$ and any power
are \Lip\ continuous on~$[\umin,\umax]$,
we see that \eqref{convforte} implies
that $f(\un)$ and $\un^{-1/2}$
converge to $f(\ui)$ and~$\ui^{-1/2}$,
respectively, in the same topology.
Therefore, we can pass to the limit 
in~\eqref{priman} 
and easily  conclude that
$\ui$ satisfies~\pbli.

\Brem
\label{Maxmonbis}
We briefly show how to extend Theorem~\ref{Longtime}
to the case mentioned in Remark~\ref{Hpmaxmon}.
Minor changes in the above proof
are necessary.
First of all, we just have to write
$f=\hat\alpha+\fd$ in~\eqref{pernuova}
and the same equality holds true.
As the sequel does not involve any
special property of~$f$,
relations \eqref{stimephi} and \eqref{stimerho}
are proved in the same way.
Hence, by going through the above proof,
we see that the only point to check
is the following:
if we define $\zetan$ by the formula
$\zetan(t)=\zeta(t+\tn)$, does 
the uniformly bounded sequence
$\graffe{\zetan}$ converge to a constant function~$\zetai$
weakly-star in~$\LQ\infty$, at least for a subsequence,
and is the pair $(\ui,\zetai)$ a solution
to the stationary problem~\pbli?
To answer these questions in the positive,
we notice that a weak-star limit $\zetai$
actually exists.
Moreover, \eqref{priman}~holds even for any $z\in\L2V$
(while so far we have written it just for 
test functions that are constant in time).
On~the other hand, the convergence~\eqref{convforte}
ensures that  $\zetai\in\alpha(\ui)$ \aeQ\
and that we can pass to the limit in~\eqref{priman}.
We obtain~that
\Beq
  \kappa \intQ \dt\ui \, z
  + \intQ \nabla\ui \cdot \nabla z
  + \intQ \bigl(\zetai + \fd'(\ui) \bigr) \, z
  - \intQ \phii \, \ui^{-1/2} \, z
  = 0
  \label{primaibis}
\Eeq
for every $z\in\L2V$.
In particular, 
$-\Delta\ui+\zetai+\fd'(\ui)=-(\phii/\ui)^{1/2}$
at least in the sense of distribution on~$Q_T$,
whence we deduce by comparison that 
$\zetai$ does not depend on time.
Then, we immediately see that 
\eqref{primaibis} becomes~\eqref{primaigen},
i.e., that $(\ui,\zetai)$ solves~\pbli.
\Erem

\section{Appendix}
\label{Sceltaxi}
\setcounter{equation}{0}
We here spend some words 
on the general forward Cauchy problem
\Beq
  y'(t) + 2a(t) \sqrt{y(t)} = 0,
  \quad
  y(0) = \yz,
  \label{ode}
\Eeq
for given $a\in\Lloc$
and $\yz\in[0,+\infty)$.

If $\yz>0$, there is a  \emph{unique, strictly positive, local solution}, 
that has the form:
\Beq
  \sqrt{y(t)} = \sqrt\yz - \iot a(s) \ds,
  \label{soluzloc}
\Eeq
as long as the \rhs\ remains positive; needless to say, $y$~may happen to tend to zero in a finite time. As to nonnegative solutions, we note that
every local solution can be extended to a global solution, because the nonlinearity is sublinear.
A~sufficient condition for uniqueness
is that $a$ is nonnegative,
because the function ~
$[0,+\infty)\ni y\mapsto a(t)\sqrt y$ ~
is nondecreasing for any fixed~$t$.
The \emph{unique, nonnegative, global solution} is given~by
\Beq
  \sqrt{y(t)}
  = \tonde{ \sqrt{\yz} - \iot a(s) \ds }^+
  \quad \hbox{for $t\in[0,+\infty),$}
  \label{soluzodesemplice}
\Eeq
where the symbol $(\cpto)^+$ denotes the positive part (for instance, if $a=1$ and $\yz=1$, 
we have that
$y(t)=((1-t)^+)^2$ for every $t\geq0$, and that $y(t)$ vanishes for $t>1$).
Actually, to verify that \eqref{soluzodesemplice}~always
provides a solution is the matter of a simple computation.
For a general~$a$, the situation is more complicated,
as we briefly explain.

If $a$ is negative, uniqueness is no longer guaranteed, because a (forward) Peano phenomenon might occur (for instance, if $a=-1$ and $\yz=0$, 
the formula $y(t)=((t-\lambda)^+)^2$
yields a solution for every $\lambda>0$).
More generally, for $\yz\geq0$,
if a solution $y^*$ vanishes a some point~$\tz$
and if $a$ is negative in a right \nbh\ of~$\tz$, then
there are infinitely many solutions beside~$y^*$.
Therefore, whenever an a priori assumption on the sign of $a$ is inappropriate (as is the case for the problem we study in our present paper),
one would like to select the solution $y$
that is \emph{maximal}, i.e., that satisfies
$y(t)\geq z(t)$ for every $t\geq0$
and for every solution $z$ to the same Cauchy problem.

Now, it is known that in general any Cauchy problem for a differential equation
of the form $y'=g(t,y)$ has a unique maximal solution
whenever $g$ is a rather general Carath\'eodory function
(see,~e.g., the extension of Peano's theorem \cite[Thm.~4.1, p.~28]{Reid}
mentioned on p.~95 of the same~book); the maximal solution can be constructed by taking
the pointwise supremum of the solutions.
For $g$ as in \eqref{ode},
the square roots of the solutions describe the~set
\Bsist
  && \calS(a,\yz) := \Bigl\{ w\in W^{1,1}_{\rm loc}[0,+\infty) :\
  w(0) = \sqrt\yz \,, 
  \enskip w(t) \geq 0 \enskip \hbox{for every $t\geq0$}, 
  \non
  \\
  && \phantom{\calS(a,\yz) := \Bigl\{} \enskip \hbox{$w'=-a$ a.e.\ where $w>0$}
  \Bigr\},
  \label{defSa}
\Esist
whence the maximal solution $y$ is \characteriz ed~by
\Beq
  \sqrt{y(t)} = \sup \graffe{w(t): \ w\in\calS(a,\yz)}
  \quad \hbox{for every $t\geq0$} .
  \label{soluzode}
\Eeq
We note that every solution 
to the Cauchy problem~\eqref{ode} is such that
\Beq
  \sqrt{y(t)} = \sqrt\yz - \iot \astar(s) \ds ,
  \label{soluzodea}
\Eeq
where
\Beq 
\astar(s) := a(s)
  \quad \hbox{if} \quad y(s) > 0
  \aand
  \astar(s) := 0 \quad \hbox{otherwise}.
  \label{soluzodeb}
\Eeq
Generally, this is nothing more than an a~posteriori reconstruction of~$y$,
because the definition of $\astar$ depends on~$y$ itself.
But, if $y$ is the maximal solution,
an additional property of $\astar$ holds true,
because $a\geq0$ a.e.\ where $y=0$.
Indeed, by maximality, $y$~tries to become positive whenever it is possible.
Thus, in view of~\eqref{soluzodeb}, we~have that
\Beq
  \astar(t) \leq a(t)
  \quad \hbox{for a.a.\ $t\geq0$}.
  \label{perpratelli}
\Eeq

At this point, we can come back to our problem. A comparison of \eqref{1pier} with~\eqref{ode} suggests that
we~take $x$ as a parameter and set:
\Beq
  a(\cpto) := \frac{\kappa \, (\dt\u(x,\cpto))^2 + \bar\emme(x,\cpto)}{2\sqrt{\u(x,\cpto)}}
  \quad \aaO\,;
  \non
\Eeq
with this,
formulas \accorpa{soluzode}{soluzodeb}
become \Seconda.
Hence, \Seconda~provide
the maximal global solution to~\eqref{1pier}
for a given~$\u$; moreover, 
a sufficient condition for uniqueness
is that $\bar\emme$ is nonnegative, a case when we have from~\eqref{soluzodesemplice} that
\Beq
  \sqrt{\xi(t)}
  = \tonde{
     \sqrt\xiz
     - \iot
        \frac{\kappa \, (\dt\u(s))^2 + \bar\emme(s)}
        {2\sqrt{\u(s)}} \ds
  }^+ .
  \label{secondasemplice}
\Eeq
If the initial datum $\xiz$ is strictly positive, the last formula holds for small $t$ even for a negative $\bar\emme$; in fact,  as long as the \rhs\ remains positive
(see  Remark~\futuro\ref{Strongsol}), we simply~have that
\Beq
  \sqrt{\xi(t)}
  = \sqrt\xiz
     - \iot
        \frac{\kappa \, (\dt\u(s))^2 + \bar\emme(s)}
        {2\sqrt{\u(s)}} \ds\,.
  \label{secondaloc}
\Eeq

\bigskip
\noindent
{\bf Acknowledgements.}\quad
Special thanks are due to Aldo Pratelli{\color{red},} who pointed out 
the proof of Lemma~\ref{ContrazODE}: the authors appreciated  
a lot his contribution. Some financial support from the MIUR-COFIN 2006 
research program on ``Free boundary problems, phase transitions 
and models of hysteresis'' and from the IMATI of CNR~in~Pavia,~Italy, 
is gratefully acknowledged. 
The work of Podio-Guidugli was supported by the Italian Ministry of
University and Research (under PRIN 2005 \textit{
``Modelli Matematici per la Scienza dei Materiali"}) 
and by the EU Marie Curie Research Training Network
MULTIMAT \emph{ ``Multi-scale Modeling and Characterization
for Phase Transformations in Advanced Materials''}.



\vspace{3truemm}

\Begin{thebibliography}{10}

\bibitem{Brezis}
H. Brezis,
``Op\'erateurs maximaux monotones et semi-groupes de contractions
dans les espaces de Hilbert'',
North-Holland Math. Stud.
{\bf 5},
North-Holland,
Amsterdam,
1973.

\bibitem{fremond}
M. Fr\'emond,
``Non-smooth Thermomechanics'',
Springer-Verlag, Berlin, 2002.

\bibitem{gurtin}
M.E. Gurtin,
{\em Generalized Ginzburg-Landau and 
Cahn-Hilliard equations based on a microforce balance},
Phys.~D
{\bf 92} (1996)
178--192.

\bibitem{Haraux}
A. Haraux,
``Syst\`emes Dynamiques Dissipatifs et Applications'',
RMA Res. Notes Appl. Math.
{\bf 17},
Masson,
Paris,
1991.

\bibitem{LSU}
O.A. Lady\v zenskaja, V.A. Solonnikov, and N.N. Ural'ceva:
``Linear and quasilinear equations of parabolic type'',
Trans. Amer. Math. Soc. {\bf 23},
Amer. Math. Soc., Providence, RI,
1968.

\bibitem{miranv} 
A. Miranville,
{\em Consistent models of {C}ahn-{H}illiard-{G}urtin equations with
{N}eumann boundary conditions},
Phys. D {\bf 158} (2001) 233--257.

\bibitem{Podio}
P. Podio-Guidugli, 
{\em Models of phase segregation and diffusion of atomic species on a lattice},
Ric. Mat. {\bf 55} (2006) 105--118.

\bibitem{Reid}
W.T. Reid,
``Ordinary differential equations'',
Wiley, New York, 1971.

\bibitem{Simon}
J. Simon,
{\em Compact sets in the space $L^p(0,T; B)$},
Ann. Mat. Pura Appl. {\bf 146} (1987) 65--96.

\End{thebibliography}

\End{document}

\bye